\documentclass[11pt,reqno]{amsart}

\usepackage{amssymb}
\usepackage{amsthm}
\usepackage{amsmath}
\usepackage[inline]{enumitem}
\PassOptionsToPackage{hyphens}{url}\usepackage{hyperref}
\usepackage{fancyhdr}
\usepackage{mathrsfs}
\usepackage{stmaryrd}
\usepackage{xcolor}
\usepackage{cite}
\theoremstyle{definition}

\newtheorem{theorem}{Theorem}[section]
\newtheorem{lemma}[theorem]{Lemma}
\newtheorem{proposition}[theorem]{Proposition}
\newtheorem{corollary}[theorem]{Corollary}

\theoremstyle{definition}

\newtheorem{example}[theorem]{Example}

\newtheorem{questions}[theorem]{Questions}

\newtheorem*{conjecture*}{Conjecture}
\newtheorem*{problem*}{Problem}
\newtheorem*{question*}{Question}

\definecolor{blue-url}{RGB}{0,0,100}
\definecolor{red-url}{RGB}{100,0,0}
\definecolor{green-url}{RGB}{0,100,0}
\definecolor{light-yellow}{RGB}{255,255,128}
\definecolor{light-blue}{RGB}{193,255,255}
\definecolor{light-red}{RGB}{239,83,80}

\hypersetup{
	pdftitle={On cancellativity and power semigroups},
	pdfauthor={Li, Tringali, and Wang},
	pdfmenubar=false,
	pdffitwindow=true,
	pdfstartview=FitH,
	colorlinks=true,
	linkcolor=blue-url,
	citecolor=green-url,
	urlcolor=red-url
}

\renewcommand{\emptyset}{\varnothing}
\renewcommand{\setminus}{\smallsetminus}
\renewcommand{\,}{\kern 0.1em}

\newcommand{\evid}[1]{\textsf{#1}}
\newcommand{\fin}{\mathrm{fin}}

\newcommand{\D}{\mathcal D}
\newcommand{\E}{\mathcal E}

\newcommand{\PP}{\mathcal P}

\title{On cancellativity and power semigroups}

\author{Lingxi Li}
\address{(L.~Li) School of Mathematical Sciences, Hebei Normal University | Shijiazhuang, Hebei province, 050024 China}
\email{lingxili.hebnu@gmail.com}

\author{Salvatore Tringali}
\address{(S.~Tringali) School of Mathematical Sciences, Hebei Normal University | Shijiazhuang, Hebei province, 050024 China}
\email{salvo.tringali@gmail.com}
\urladdr{http://imsc.uni-graz.at/tringali}

\author{Chen Wang}
\address{(C.~Wang) School of Mathematical Sciences, Hebei Normal University | Shijiazhuang, Hebei province, 050024 China}
\email{chen.wang.hebnu@outlook.com}

\date{}

\begin{document}
\begin{abstract}
The large power semigroup of a semigroup $S$ is the semigroup of all non-empty subsets of $S$, equipped with the setwise extension of its operation.

We prove that every isomorphism between the large power semigroups of two cancellative semigroups $H$ and $K$ maps singletons to singletons and hence induces an isomorphism from $H$ to $K$. As a result, we settle a question posed by Garc\'ia-S\'anchez and Tringali in~[Proc.\ Amer.\ Math.\
Soc.\ \textbf{153} (2025), 2323--2339].

The proof is essentially combinatorial: it exploits cancellativity to characterize singletons in terms of the cardinalities of
the solution sets of certain equations of the form $AX = B$. The same argument applies more generally to isomorphisms between subsemigroups
of the large power semigroups of $H$ and $K$, on condition that these subsemigroups are closed under taking binary unions and non-empty subsets (we
provide several examples).
\end{abstract}

\maketitle

\section{Introduction}\label{sec:introduction}

Let $S$ be a semigroup (see the end of the section for terminology and notation). As usual, we define the setwise product of two sets $X, Y \subseteq S$ by
\[
(X, Y) \mapsto XY := \{xy: x \in X,\ y \in Y\},
\]
Under this operation, the non-empty subsets of $S$ form a semigroup $\PP(S)$, called the \evid{large power semigroup} of $S$. The subsemigroup of $\mathcal P(S)$ consisting of all non-empty finite subsets of $S$ is the \evid{finitary power semigroup} of $S$ and is denoted by $\PP_\fin(S)$.

The systematic investigation of power semigroups began with Tamura and Shafer~\cite{Tam-Sha-1967} in the late 1960s. After a period of intense activity from the 1970s through the 1990s, interest in the topic waned. In recent years, however, research on the subject has regained momentum, with a growing body of work addressing a broad range of
questions, notably including the study of automorphism groups~\cite{Tri-Yan2025(b), Tri-Wen-2026(a), Rago26(b), Tri-Wen-2026(b), Wong2026}, isomorphism problems~\cite{Bie-Ger-22, Tri-Yan2025(a), GarSan-Tri-24(a), Rago26(a),
Rago26(c), Tri-Yan2026(a)}, and arithmetic properties~\cite{Fa-Tr18, Ant-Tri-PJM2018, Bie-Ger-22, Co-Tr-26(a)}. We refer the reader to \cite{Tri-AMM2027} for a survey of these recent developments.

Tamura, in particular, maintained a lifelong interest in the first
of the following questions (here and throughout, $\cong$ denotes
semigroup isomorphism).

\begin{questions}
\label{que:tamura-shafer}
Let $\mathcal O$ be a class of semigroups. Prove or disprove that, for all $H, K \in \mathcal O$, 
\begin{enumerate}[label=\textup{(\arabic{*})}]
\item\label{que:tamura-shafer(1)}
$\mathcal P(H) \cong \mathcal P(K)$ if and only if $H \cong K$.
\item\label{que:tamura-shafer(2)} $\mathcal P_\fin(H) \cong \mathcal P_\fin(K)$ if and only if $H \cong K$.
\end{enumerate}
\end{questions}

The converse implication is trivial in either case. Indeed, an isomorphism $\phi \colon H \to K$ induces an isomorphism $\phi^\ast \colon \mathcal P(H) \to \mathcal P(K)$ that sends a non-empty $X \subseteq H$ to its direct image $f[X] := \allowbreak \{f(x): x \in X\}$, and restricting $\phi^\ast$ to $\mathcal P_\fin(H)$ yields an isomorphism onto $\mathcal P_\fin(K)$. 

A class $\mathcal{O}$ of semigroups is \evid{globally determined} if Question~\ref{que:tamura-shafer}\ref{que:tamura-shafer(1)} has a positive answer for $\mathcal{O}$. The term derives from the fact that, in the older literature, large power semigroups were also called \evid{globals}. In keeping with this terminology, for semigroups $H$ and $K$, an isomorphism from $\mathcal P(H)$ to $\mathcal P(K)$ is still sometimes referred to as a \evid{global isomorphism} from $H$ to $K$. 

While the class of all semigroups is not globally determined~\cite{Mog-1973}, several important subclasses are, including groups~\cite{Sha-1967}, semilattices~\cite[p.~218]{Kob-1984}, and, more generally, Clifford semigroups~\cite[Theorem~4.7]{Gan-Zha-2014}. In each of these cases, the proof establishes a stronger conclusion: every global isomorphism maps singletons to singletons and hence induces an isomorphism between the underlying semigroups. This is sometimes referred to as the \evid{Kobayashi property} and leads to the following variant of Questions~\ref{que:tamura-shafer}.

\begin{questions}
\label{que:kobayashi}
Let $\mathcal O$ be a class of semigroups. Prove or disprove that, for all $H, K \in \mathcal O$, 
\begin{enumerate}[label=\textup{(\arabic{*})}]
\item\label{que:kobayashi(1)}
every isomorphism $\mathcal P(H) \to \mathcal P(K)$ preserves singletons.
\item\label{que:kobayashi(2)} every isomorphism $\mathcal P_\fin(H) \to \mathcal P_\fin(K)$ preserves singletons.
\end{enumerate}
\end{questions}

Both item~\ref{que:kobayashi(1)} and \ref{que:kobayashi(2)} of Questions~\ref{que:kobayashi} have positive answers for inverse semigroups~\cite[Corollaries~4.7 and 5.4]{Liu-Tri-2026(b)} and for cancellative \emph{commutative} semigroups~\cite[Corollary~1]{Tri-2024(a)}. In the present paper, we generalize the latter result by showing that the answer is still affirmative for arbitrary cancellative semigroups (Corollary~\ref{cor:powers}). This settles an open question posed by Garc\'ia-S\'anchez and the second author in~\cite[p.~2325, ll.~3--4]{GarSan-Tri-24(a)}.

The proof is combinatorial in nature: it characterizes singletons in terms of the cardinalities of the solution sets of certain equations of the form $AX=B$, exploiting the cancellativity of the
underlying semigroups. This makes it possible to overcome one of the main difficulties in the Tamura--Shafer problem and its many variants: isomorphisms between power semigroups need not be inclusion- or cardinality-preserving. 

The same argument applies more generally to isomorphisms between hereditary, union-closed subsemigroups
of the large power semigroups of cancellative semigroups (Theorem~\ref{thm:singletons}). Here, a family $\mathcal D$ of sets is called \evid{union-closed} if $X \cup Y \in \D$ for all $X, Y \in \D$, and \evid{hereditary} if it contains all non-empty subsets of any set in the family. Notice that if $\D$ is hereditary, then it contains each one-element subset of $\bigcup \D := \bigcup_{X \in \D} X$. In particular, $\mathcal P(S)$ and $\mathcal P_\fin(S)$ are hereditary and union-closed for every semigroup $S$. Further examples are discussed below.

\begin{example}
\label{example:01}
Let $\kappa$ be an infinite cardinal, and consider the families
\[
\PP_{\le \kappa}(S) := \{X \in \PP(S): |X| \le \kappa\}
\quad\text{and}\quad
\PP_{<\kappa}(S) := \{X \in  \PP(S):|X| < \kappa\}. 
\]
Both are hereditary,
union-closed subsemigroups of $\PP(S)$. Indeed, the cardinality bounds are preserved under taking subsets and finite unions. Moreover,
$|XY|\le |X|\,|Y|$ for all $X, Y \subseteq S$.
\end{example}

\begin{example}\label{example:02}
Suppose that $S$ is a semigroup, and let
$\ell$ be a function $S \to \mathbb{R}_{\ge 0}$ such that
\[
\ell(xy) \le \ell(x) + \ell(y),
\qquad\text{for all }
x, y \in S. 
\]
That is, $\ell$ satisfies the \evid{triangle inequality}. For every $X \subseteq S$, set
\[
m_\ell(X) := \sup_{x \in X} \ell(x)
\quad\text{and}\quad
\mathcal{B}_\ell(S) := \{X \in \PP(S): m_\ell(X) < \infty\},
\]
where $\sup \emptyset := 0$. If $X, Y \subseteq S$, then $X \subseteq Y$ implies $m_\ell(X) \le m_\ell(Y)$. In addition, we have
\[
m_\ell(X \cup Y) \le \max (m_\ell(X), m_\ell(Y))
\quad\text{and}\quad
m_\ell(XY) \le m_\ell(X) + m_\ell(Y).
\]
It follows that $\mathcal{B}_\ell(S)$ is a hereditary, union-closed subsemigroup of $\PP(S)$. 

If, in particular, $S$ is the additive group of the real field and $\ell$ is the usual absolute
value on $\mathbb{R}$, then $\mathcal{B}_\ell(S)$ consists of the non-empty bounded subsets of $\mathbb{R}$. Unlike the constructions in Example~\ref{example:01}, $\mathcal{B}_\ell(S)$ is not determined by cardinality bounds. Indeed, $[0,1]$ belongs to $\mathcal{B}_\ell(S)$,
whereas $\mathbb{R}$ does not, although the two sets have the same
cardinality.
\end{example}

The semigroups in Example~\ref{example:01} were already considered in \cite[Example~1]{Tri-2024(a)}, whereas those in Example~\ref{example:02} appear to be new. It may be worthwhile to study them in greater depth.

\subsection*{Generalities.} We denote by $\mathbb N$ the (set of) non-negative integers, by $\mathbb N^+$ the positive integers, by $\mathbb R$ the real numbers, and by $|X|$ the cardinality of a set $X$.

We refer to~\cite{cli-pre-1961-I} for basic aspects of semigroup theory. All morphisms are understood to be semigroup homomorphisms. Moreover, by a submonoid (resp., subgroup) of a semigroup $S$ we mean a subsemigroup $T$ of $S$ that is a monoid (resp., a group) in its own right. In particular, if $S$ is itself a monoid, we do not require the identity of $T$ to coincide with that of $S$.  When they do coincide, we say that $T$ is a \evid{unitary} submonoid (resp., subgroup) of $S$.

If not explicitly specified, we write all semigroups multiplicatively. An element $a$ in a semigroup $S$ is \evid{cancellative} if $ax \ne ay$ and $xa \ne ya$ for all $x, y \in S$ with $x \ne y$; the semigroup itself is cancellative if each of its elements is. The \evid{order} of an element $a \in S$ is the cardinality of the cyclic subsemigroup $\{a^n: n \in \mathbb N^+\}$ generated by $a$. A \evid{unit} in a monoid $T$ is an element $u \in T$ for which there exists a (provably unique) element $v \in T$ such that $uv = vu = 1_T$. The units form a subgroup of $T$, which we denote by $T^\times$ and call the \evid{unit group} of $T$. Since $1_T \in T^\times$, it is clear that $T^\times$ is a unitary subgroup of $T$.

Any further notation or terminology, if not explained when first introduced, is standard or should be clear from the context. 

\section{Counting solutions}\label{sec:counting}

Our goal in this section is to obtain an \textit{algebraic} characterization of the singletons in hereditary, union-closed subsemigroups of large power semigroups of cancellative semigroups
(Proposition~\ref{prop:singletons}). To this end, given a semigroup $S$ and elements $a,b\in S$, we define
\[
N_S(a,b) := \bigl|\{x \in S: ax = b\}\bigr|.
\]
The insight here is that isomorphisms preserve the cardinality of the solution set to every equation of the form $ax = b$. The idea is not new, but was never used before in conjunction with cancellativity. We begin with two lemmas.

\begin{lemma}\label{lem:deletion}
Let $S$ be a cancellative semigroup and $\D$ be a hereditary, union-closed subsemigroup of $\PP(S)$. Fix $A, B, C \in \mathcal D$ with $|C| \ge 2$. If $N_\D(A, B) \ge 2$, then
$N_{\D}(A,BC) \ge 3$.
\end{lemma}

\begin{proof}
Assume $N_\D(A, B) \ge 2$. The solutions in $\D$ of the equation $AZ = B$ are closed under binary unions. Indeed, $\D$ is union-closed. Thus, if $AX = AY = B$ for some $X, Y \in \D$, then
\[
X \cup Y \in \D
\quad\text{and}\quad
A(X \cup Y) = AX \cup AY = B.
\]

Choose two distinct solutions $X$ and $Y$. Their union properly contains at least one of them, so we may assume without loss of generality that
\[
X \subsetneq Y
\quad\text{and}\quad
AX = AY = B.
\]
Pick $y \in Y \setminus X \ne \emptyset$, and let $c \in C$. Since $\D$ is a subsemigroup, it is clear that
\begin{equation}\label{lem:deletion:eq(1)}
Q := YC \in \D
\quad\text{and}\quad
AQ = (AY)C = BC.
\end{equation}
Also, right cancellativity yields $yc \notin Xc$, and hence
\begin{equation}\label{lem:deletion:eq(2)}
\emptyset \ne Xc \subseteq (YC) \setminus \{yc\} = Q \setminus \{yc\} =: Q_c.
\end{equation}
Accordingly, we find that
\begin{equation}\label{lem:deletion:eq(3)}
Ayc = (Ay)c \subseteq (AY)c = (AX)c
= A(Xc) \subseteq AQ_c.
\end{equation}
Upon considering that $Q = Q_c \cup \{yc\}$, it follows by Eqs.~\eqref{lem:deletion:eq(1)} and \eqref{lem:deletion:eq(3)} that
\[
BC = AQ = Ayc \cup AQ_c = AQ_c.
\]
In addition, Eq.~\eqref{lem:deletion:eq(2)} implies that $Q_c$ is non-empty, so $Q_c \in \D$ by the hereditarity of $\mathcal D$.

Now choose distinct elements $c_1, c_2 \in C$, which is possible since $|C| \ge 2$. By left cancellativity, $yc_1 \ne yc_2$, and hence
$Q_{c_1} \ne Q_{c_2}$. Since each $Q_{c_i}$ is properly contained in $Q$, the sets $Q$, $Q_{c_1}$, and $Q_{c_2}$ are therefore pairwise distinct solutions in $\D$ of
$AZ = BC$. So, $N_{\D}(A,BC)\ge 3$.
\end{proof}

\begin{lemma}\label{lem:two-solutions}
Let $S$ be a cancellative semigroup, let $a \in S$ be an element of
infinite order, and let $b \in S$. For every integer $n \geq 1$, the
equation
\begin{equation}\label{lem:two-solutions:eq(1)}
\{a, a^2\}X = \{a, \ldots, a^n\}b
\end{equation}
has exactly $F_{n-1}$ solutions $X \in \mathcal{P}(S)$, where $F_0, F_1, \ldots$ is the Fibonacci sequence (with 
$F_0 := 0$ and $F_1 := 1$). Moreover, any such solution is finite. 
\end{lemma}

\begin{proof}
Let $\widehat{S}$ be the monoid obtained from $S$ by adjoining an identity element only if $S$ is not already a monoid, and let $e$ be the identity element of $\widehat{S}$. Since $S$ is cancellative, so is $\widehat{S}$. 

Building on these preliminaries, let $X \in \mathcal P(S)$, and fix an integer $n \geq 1$. Factoring out $\{a\}$ on the left and using the left cancellativity of $\widehat{S}$, we find that
$$
\{a, a^2\}X = \{a, \ldots, a^n\}b,\, \text{ in } \mathcal P(S)
\quad\Longleftrightarrow\quad
\{e, a\}X = \{e,a\}^{n-1} b,\, \text{ in } \mathcal P(\widehat{S}).
$$
Thus, $X$ is a solution of Eq.~\eqref{lem:two-solutions:eq(1)} only if
\begin{equation}\label{lem:two-solutions:eq(2)}
X = eX \subseteq \{e,a\}X = \{e,a\}^{n-1} b.
\end{equation}
In particular, Eq.~\eqref{lem:two-solutions:eq(1)} has finitely many solutions, and each of them is finite. Moreover, right cancellativity of $\widehat{S}$ gives a unique non-empty subset $Y \subseteq \{e,a\}^{n-1}$ such that $X = Yb$. Substituting
$Yb$ and cancelling $b$ on the right, we obtain
\[
\{e,a\}X =\{e,a\}^{n-1}  b,\, \text{ in } \mathcal P(\widehat{S})
\quad\Longleftrightarrow\quad
\{e,a\}Y = \{e,a\}^{n-1},\, \text{ in } \mathcal P(\widehat{S}).
\]
Conversely, it is clear that every solution $Y \in \mathcal{P}(\widehat{S})$ of the
last equation gives a solution $Yb \in \mathcal{P}(S)$ of Eq.~\eqref{lem:two-solutions:eq(1)}. In other words, the two equations have the same number $s_n$ of
solutions. Since $a$ has infinite order in $\widehat{S}$, we therefore conclude from \cite[Proposition 3.2]{Liu-Tri-2026(a)} that $s_n = F_{n-1}$. 
\end{proof}

The following proposition is the key ingredient in the proof of our main result (Theorem \ref{thm:singletons}).

\begin{proposition}\label{prop:singletons}
Let $S$ be a cancellative semigroup and $\D$ be a hereditary,
union-closed subsemigroup of $\PP(S)$. Assume that $S$ has at least one element of infinite order. A set $Y \in \D$ is then a singleton
if and only if there exists $A \in \D$ such that
\begin{equation}
\label{prop:singletons:eq(1)}
 N_{\D}(A, A^3) = N_{\D}(A, A^3 Y) = 2.
\end{equation}
\end{proposition}

\begin{proof}
Suppose that $Y$ is a singleton, say, $Y = \{y\}$. By hypothesis, there is at least one element $a \in S$ of infinite order. Set $A := \{a,a^2\}$ and $B := \{a,a^2,a^3,a^4\}$. We have
\[
 A^3 = \{a^3, a^4, a^5, a^6\} = Ba^2
 \quad\text{and}\quad
 A^3 Y = Ba^2y.
\]
Now apply Lemma~\ref{lem:two-solutions}, first with $b = a^2$ and then with $b = a^2y$. In each case, the equation $AX = Bb$ has two solutions in $\mathcal P(S)$, and both of them are finite. Since $\mathcal D$ is hereditary and hence contains every non-empty subset of its own elements, we conclude that 
\[
N_{\D}(A, A^3) = N_{\D}(A, A^3Y) = 2.
\]

Conversely, suppose there exists $A \in \D$ such that Eq.~\eqref{prop:singletons:eq(1)} holds.
If $|Y| \ge 2$, then applying Lemma~\ref{lem:deletion} with $B = A^3$ gives $N_{\D}(A, A^3Y) \ge 3$, which is absurd. Thus, $Y$ is a singleton.
\end{proof}

\section{Main result}\label{sec:main-result}

We are ready to prove that both parts of Questions~\ref{que:kobayashi} have positive answers for the class of cancellative semigroups. We begin with a proposition
that generalizes Remark~2.2(2) in~\cite{Liu-Tri-2026(a)} and
Theorem~4.2(ii) in~\cite{Tri-AMM2027}. 

\begin{proposition}\label{prop:monoid-vs-monoid}
Let $S$ be a semigroup, and let $\D$ be a hereditary submonoid of
$\PP(S)$ such that $\bigcup\D = S$. Then $S$ is a monoid, with $1_{\D} = \{1_S\}$. Moreover, the units of $\D$ are precisely the one-element subsets of the set of units of $S$.
\end{proposition}

\begin{proof}
Let $E$ be the identity of $\D$, and choose $e \in E$ (this is possible because every set in $\D$ is non-empty). Since $\D$ is hereditary and $S$ is the union of the sets in $\D$, we have that
\begin{equation}\label{prop:monoid-vs-monoid:eq(1)}
\{x\} \in \D,
\qquad\text{for every }x \in S.
\end{equation}

For each $x \in S$, we have $Ex = xE = \{x\}$, and hence $ex = xe = x$. Therefore, $e$ is an identity of $S$. Moreover, the same argument shows that any element of $E$ is indeed an
identity of $S$. By uniqueness of the identity, it follows that $E = \{e\}$. Consequently, $S$ is a monoid with identity $1_S = e$. In addition, we obtain $1_{\D} = E = \{1_S\}$.

Denote by $S^\times$ the set of units of $S$ and recall from Eq.~\eqref{prop:monoid-vs-monoid:eq(1)} that $\D$ contains every one-element subset of $S$. If $u \in S^\times$, then $\{u\}$ is a unit of $\D$, its inverse being the singleton $\{u^{-1}\}$.
Conversely, let $X$ be a unit of $\D$, and let $Y$ be an inverse of $X$ in $\D$, so that
$XY = YX = \{1_S\}$. Fix $y \in Y$. For every $x \in X$, we have
$xy = yx = 1_S$, which shows that $y$ is a unit of $S$ and $x$ is an inverse of $y$. By uniqueness of the inverse, this proves that $X$ is a singleton, as required.
\end{proof}

We are finally ready to prove the main results of the paper.

\begin{theorem}\label{thm:singletons}
Let $H$ and $K$ be cancellative semigroups, and let
$\D$ and $\E$ be hereditary,
union-closed subsemigroups of $\PP(H)$ and $\PP(K)$, respectively.
Every isomorphism $\Phi \colon \D \to \E$ maps singletons to singletons.
\end{theorem}

\begin{proof}
It is clear that $H' := \bigcup \mathcal D$ and $K' := \bigcup \mathcal E$ are subsemigroups of $H$ and $K$, respectively. In particular, they are both cancellative. Since $\D \subseteq \PP(H')$ and $\E \subseteq \PP(K')$, we may therefore assume without loss of generality that $\emptyset \ne H = H'$ and $\emptyset \ne K = K'$ (if $H$ or $K$ is empty, then so are $\D$ and $\E$, which makes the conclusion trivial). 

Suppose first that every element of $H$ has finite order. It is a basic fact that $H$ is then a group; see, for instance, Exercise~6(c) in \cite[Section~1.7]{cli-pre-1961-I}. By the hypothesis that $\D$ is hereditary (and $H$ is union of the sets in $\D$), every one-element subset of $H$ belongs to $\D$. In particular, $\D$ is a submonoid of $\PP(H)$ with identity $\{1_H\}$. So, by Proposition~\ref{prop:monoid-vs-monoid}, the units of $\D$ are precisely the one-element subsets of $H$. Since $\Phi$ is an isomorphism $\D \to \E$, it follows that $\E$ is a monoid with identity
$\Phi(\{1_H\})$, and $\Phi$ maps $\D^\times$ onto $\E^\times$. This proves that $\Phi$ maps singletons to singletons, because Proposition~\ref{prop:monoid-vs-monoid} applied to $\E$ shows that $K$ is itself a monoid and the units of $\E$ are precisely the one-element subsets of $K^\times$.

Now assume that $H$ has at least one element of infinite order, say $a$. Given $C \in \D$, consider the condition that there exists
$A \in \D$ such that
\begin{equation*}
N_{\D}(A, A^3) = N_{\D}(A, A^3C) = 2.
\end{equation*}
For any $A, B \in \D$, the isomorphism $\Phi$ maps the solutions in $\D$ of $AX = B$ bijectively onto the solutions in $\E$ of $\Phi(A)Y = \Phi(B)$. Since $\Phi(A^3C) = \Phi(A)^3\Phi(C)$, it follows that $\Phi$ preserves the condition in the last display. On the other hand, we have from Proposition~\ref{prop:singletons} that the same condition characterizes the
singletons in $\D$. In particular, $\{a\}\in\D$ satisfies it, so
$\Phi(\{a\})$ does the same in $\E$. The implication from the condition to being
a singleton follows from Lemma~\ref{lem:deletion}, even without assuming
that $K$ contains an element of infinite order. 

Consequently, $\Phi(\{a\}) = \{b\}$ for some $b\in K$. Since an isomorphism preserves the order of the elements, $b$ has infinite order in $K$. Thus $K$ also contains an element of infinite
order. We may therefore apply Proposition~\ref{prop:singletons} to both $H$ and
$K$. The proposition identifies the singletons in $\D$ and $\E$ by the condition above, which is preserved by $\Phi$. Therefore, $\Phi$ maps singletons to singletons, as required.
\end{proof}

\begin{corollary}\label{cor:powers}
If $H$ and $K$ are cancellative semigroups, then every isomorphism $\mathcal P(H) \to \mathcal P(K)$ or $\mathcal P_\fin(H) \to \mathcal P_\fin(K)$ maps singletons bijectively onto singletons. In particular, $H \cong K$.
\end{corollary}

\begin{proof}
For every semigroup $S$, both $\PP(S)$ and $\PP_\fin(S)$ are hereditary and union-closed, and each of them contains every one-element subset of $S$. Since the inverse of an isomorphism is itself an isomorphism, the conclusion is then immediate from Theorem~\ref{thm:singletons}.
\end{proof}

\section*{AI disclaimer}

The idea of characterizing the one-element subsets of the large power semigroup of a cancellative semigroup by means of the condition in Eq.~\eqref{prop:singletons:eq(1)} was suggested in a conversation with AI, possibly prompted by Proposition~3.2 of~\cite{Liu-Tri-2026(a)}. The underlying strategy of counting solutions to equations of the form $AX = B$ is not new; Tringali et al.~had already used it in earlier work. There is, however, some ingenuity in the way AI adapted this strategy to the cancellative setting. The authors are responsible for the subsequent development of this idea, including its extension to hereditary, union-closed subsemigroups, and for all other contributions in this paper.
\section*{Acknowledgments}

The authors were supported by the 111 Center under Grant No.~D26018 and by the Natural Science Foundation of Hebei Province under Grant No.~A2023205045.


\end{document}